\documentclass[11pt]{article}
\usepackage{amsmath}
\usepackage{amssymb}
\usepackage{amsthm}
\usepackage{amsfonts}
\usepackage{pst-plot}

\newcommand{\piece}[1]{\bigskip \noindent \emph{#1}}
\newcommand{\step}[1]{\bigskip \emph{Step #1.}}
\newcommand{\case}[1]{\bigskip \emph{Case #1.}}

\newenvironment{mainthm}{\bigskip \noindent {\bf Main Theorem. }\em{}{}}

\DeclareMathOperator{\mcg}{Mod}
\DeclareMathOperator{\aut}{Aut}
\DeclareMathOperator{\isom}{Isom}
\DeclareMathOperator{\out}{Out}
\DeclareMathOperator{\teich}{\mathcal{T}}

\DeclareMathOperator{\cpo}{C_P^1}
\DeclareMathOperator{\cpos}{C_P^1(S)}
\DeclareMathOperator{\cps}{C_P(S)}
\DeclareMathOperator{\cc}{C}
\DeclareMathOperator{\ccz}{C^{(0)}}

\newcommand{\al}{\ensuremath{\alpha}}
\newcommand{\be}{\ensuremath{\beta}}
\newcommand{\ga}{\ensuremath{\gamma}}

\newcommand{\ep}{\ensuremath{\epsilon}}

\newcommand{\lp}{\ensuremath{\mathcal{L}}}

\newcommand{\szfi}{\ensuremath{\Sigma_{0,5}}}
\newcommand{\szfo}{\ensuremath{\Sigma_{0,4}}}
\newcommand{\sot}{\ensuremath{\Sigma_{1,2}}}
\newcommand{\soo}{\ensuremath{\Sigma_{1,1}}}

\newcommand{\mods}{\ensuremath{\mcg(S)}}

\newcommand{\lra}{\ensuremath{\longrightarrow}}
\newcommand{\ra}{\ensuremath{\rightarrow}}
  
\newcommand{\twps}{\ensuremath{\teich_{WP}(S)}}  
\newcommand{\tts}{\ensuremath{\teich (S)}}  
\newcommand{\ttss}{\ensuremath{\teich (S)\ }}  
\newcommand{\twpss}{\ensuremath{\teich_{WP}(S)\ }}  
\newcommand{\tm}{Teichm\"{u}ller\,}

\newcommand{\tss}{\tm Space\ }

\begin{document}

\input{epsf.sty}

\title{\bf The Automorphism Group of the Complex of Pants Decompositions}
\author{Dan Margalit}

\maketitle

\section{Introduction}

There has been recent interest in finding combinatorial models for groups, i.e. finding a simplicial complex whose
automorphism group is a given group.  For example, Charney-Davis showed that the Coxeter Diagram is a model for
$\out(W)$ (for certain Coxeter groups W) \cite{cd}, and Bridson-Vogtmann showed that a spine of Outer Space is a
model for $\out(F_n)$ \cite{bv}.  In this paper, we prove that the complex of pants decompositions is a model for the
mapping class group.

Throughout, $S$ will denote a closed oriented surface of negative Euler characteristic.

The complex of pants decompositions of S, defined by Hatcher-Thurston and denoted $\cps$, has vertices
representing pants decompositions of $S$, edges connecting vertices whose pants decompositions differ by an 
elementary move (see below), and 2-cells representing certain relations between elementary moves.  Its 1-skeleton is
called the graph of pants decompositions, and is denoted $\cpos$.

Brock proved that $\cpos$ models the \tss of $S$ endowed with the Weil-Petersson metric, $\twps$, in that the spaces
are quasi-isometric \cite{jb}.  Our main result further indicates that $\cpos$ is the ``right'' combinatorial model
for \twps, in that the automorphism group of $\cpos$ is shown to be the extended mapping class group of $S$, $\mods$
(the group of diffeomorphisms of $S$ to itself, modulo isotopy).  This is in consonance with the result of Masur-Wolf
that the isometry group of \twpss is $\mods$ \cite{mw}.  The extended mapping class group has a natural action 
by automorphisms on $\cpos$; the content of the theorem is that all of the automorphisms of $\cpos$ are induced by
elements of $\mods$.

The Main Theorem of this paper is as follows:

\begin{mainthm}Let $S$ be a closed oriented surface of negative Euler Characteristic which is not the genus $1$
surface with $2$ punctures.  Then:

\[ \aut \cps \cong \aut \cpos \cong \mods \]

\end{mainthm}

The two isomorphisms are proved in Sections ~\ref{graph-complex} and ~\ref{graph-curvecomplex}, respectively.  To
prove the first isomorphism of the Main Theorem, we show that the 2-cells of $\cps$, which are defined via
topological relations between pants decompositions on $S$, can equivalently be characterized using only the
combinatorics of $\cpos$.  In other words, $\cps$ carries no more information than its 1-skeleton. For example,
square 2-cells of $\cps$ (those representing relations involving four elementary moves) are defined as a commutator
relation made up of elementary moves on disjoint subsurfaces on $S$ (see Figure~\ref{squpic}). We show that square
2-cells can equivalently be defined as loops with 4 edges in $\cpos$ which have the property that adjacent edges in
the loop do not lie in a common Farey graph (see below) in $\cpos$. Note that in the second definition there is no
mention of the surface, only $\cpos$.  Therefore, any automorphism of $\cpos$ preserves these square 2-cells.

In order to prove the second isomorphism of the Main Theorem, we use the fact that for any $S$ satisfying the
hypotheses of the Main Theorem, we have the following isomorphism:

\[ \aut \cc (S) \cong \mods \]

\noindent where $\cc(S)$ is the complex of curves for $S$ (defined below).  This isomorphism is a theorem of Ivanov
\cite{iv}. Korkmaz proved the low genus case \cite{ko}, and Luo gave another proof for any genus \cite{fl}.

While $\cpos$ has a vertex for each top-dimensional simplex of $\cc(S)$, it is only a relatively small subcomplex of
the dual of $\cc(S)$, and hence the main result of this paper does not follow from Ivanov's theorem.

The key idea for the second isomorphism of the Main Theorem is that there is a correlation between marked Farey
graphs (see below) in $\cpos$ and vertices in $\cc(S)$.  An automorphism of $\cpos$ induces a permutation of these
Farey graphs, and hence gives rise to an automorphism of $\cc(S)$.  The difficulty is to show that this map between
the automorphism groups of $\cpos$ and $\cc(S)$ is well-defined.

\section{The Complexes} \label{complexes}

\paragraph{Complex of Curves.} The {\em complex of curves} (or {\em curve complex}) of $S$ is the simplicial complex
$\cc(S)$ whose vertices correspond to isotopy classes of nontrivial simple closed curves on $S$;  a curve is
nontrivial if it is essential (not null homotopic) and nonperipheral (not homotopic to a boundary component).  
Throughout, we will use {\em curve} to mean {\em isotopy class of curves}.  Also we will lose the distinction between
isotopy classes of curves and representatives of the isotopy classes.

A set of $k+1$ vertices of the curve complex is the $0$-skeleton of a $k$-simplex if the corresponding curves have
trivial intersection pairwise (there is a set of representatives of the classes which are mutually disjoint).  For
example, edges correspond to pairs of disjoint curves, triangles correspond to triples of disjoint curves, etc.

The curve complex was first defined by Harvey \cite{wh}.  Harer proved that it is homotopy equivalent to a wedge of
spheres \cite{jh}.  Ivanov used his theorem that $\aut \cc (S) \cong \mods$ to give a new proof of Royden's Theorem
that $\isom(\tts) \cong \mods$ (where \ttss is the \tss of $S$ with the \tm metric) \cite{iv}.  Masur-Minksy have
shown that the complex of curves is $\delta$-hyperbolic \cite{mm}.

\begin{figure}[htb]
\begin{center}
\input{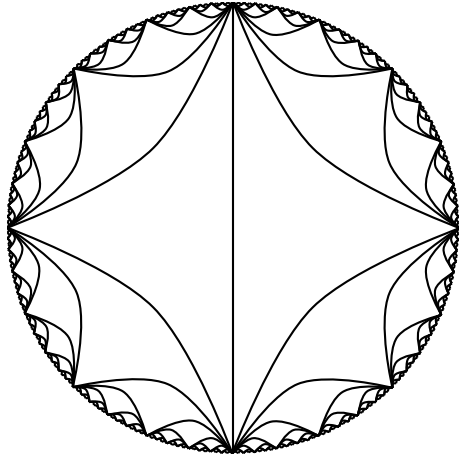}
\end{center}
\caption{An ideal triangulation of the disk, or {\em Farey graph}.
\label{farey}
}
\end{figure}

The complex of curves has an altered definition in two cases.  For $\szfo$ and $\soo$ ($\Sigma_{g,b}$ denotes the
genus $g$ surface with $b$ boundary components), since there are no distinct simple closed curves with trivial
intersection, two vertices are connected by an edge when the curves they represent have minimal intersection (2 in
the case of $\szfo$, and 1 in the case of $\soo$).  It turns out that in both cases, the complex of curves is an
ideal triangulation of the disk, or {\em Farey Graph} (see Figure~\ref{farey}) \cite{mi}.

\paragraph{Pants Decompositions.} A {\em pants decomposition} of $S$ is a maximal collection of distinct nontrivial
simple closed curves on $S$ which have trivial intersection pairwise. In other words, pants decompositions correspond
to maximal simplices of the curve complex.  A pants decomposition always consists of $3g-3+r$ curves (where $S$ is
genus $g$ with $r$ punctures). The complement in $S$ of the curves of a pants decomposition is $2g-2+r$ thrice
punctured spheres, or {\em pairs of pants}.  A pants decomposition is written as $\{\al_1, \dots,\al_n\}$, where the
$\al_i$ are curves on $S$.

Two curves in a pants decomposition are said to {\em lie on disjoint subsurfaces} if they are not the boundary
components of a common pair of pants in $S$.

\paragraph{Elementary Moves.} Two pants decompositions $p$ and $p'$ of $S$ are said to {\em differ by an elementary
move} if $p'$ can be obtained from $p$ by replacing one curve, say $\al_1$ in $p$, with another curve, say $\al_1'$,
such that $\al_1$ and $\al_1'$ intersect \emph{minimally}---if $\al_1$ lies on a $\szfo$ in the complement of the
other curves in $p$, then {\em minimally} means $\al_1$ and $\al_1'$ must intersect exactly twice; if $\al_1$ lies on
a $\sot$ in the complement of the rest of $p$, $\al_ 1$ and $\al_1'$ should intersect exactly once.  These are the
only possibilities, corresponding to whether $\al_1$ is the boundary between two pairs of pants on $S$ or is on a
single pair of pants.

In the case of $\szfo$ the curve involved in the elementary move changes its {\em
association}: the association of a curve on $\szfo$ is the natural grouping it induces on the
punctures---the curve associates two punctures if they lie on the same subsurface in the complement of the curve.  We
note that given a curve $\al$ on $\szfo$ with a certain association and some different specified association,
then up to Dehn twists (see \cite{iv2}) about $\al$ and a reflection which fixes $\al$, there is a unique choice of
curve $\be$ which has the specified association and which differs from $\al$ by an elementary move.  We call
this the {\em associativity move rule}, and will refer to it several times in the proof of the Main Theorem.

An elementary move will be denoted $p \ra p'$, $\{\al_1,\dots,\al_n\} \ra \{\al_1',\al_2,\dots,\al_n\}$, or $\al_1
\ra \al_1'$.  Note that there are countably many elementary moves $\al_1 \ra \star$. 

\epsfysize=2 cm
\begin{figure}[htb]
\center{
\leavevmode
\epsfbox{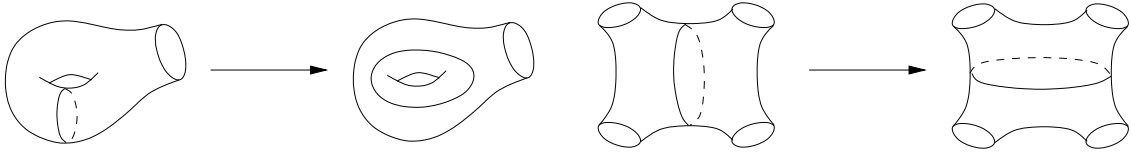}
\caption{Elementary moves between pants decompositions.}
}
\end{figure}

\paragraph{Graph of Pants Decomositions.} The {\em graph of pants decompositions} (or {\em pants graph}) of $S$,
denoted $\cpos$, is the abstract simplicial complex with a vertices corresponding to pants
decompositions of $S$, and edges joining vertices whose associated pants decompositions differ by an
elementary move.

Note that the pants graphs for $\szfo$ and $\sot$ have the same definitions as (the 1-skeleton of) the curve
complexes for these surfaces---all four complexes are Farey graphs.

\begin{figure}[htb]
\center{
\leavevmode
\epsfbox{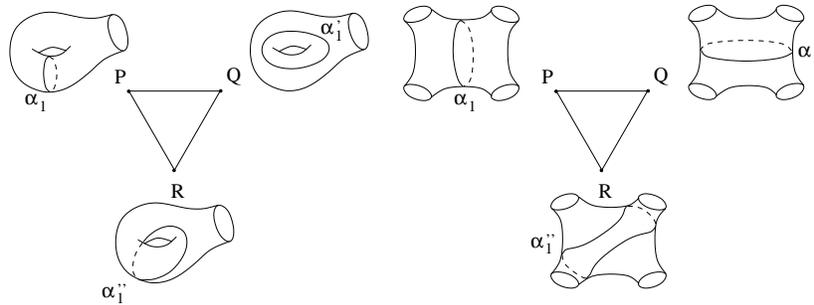} 
\caption{Triangles in the Pants Complex.}
\label{tripic}
}
\end{figure}

\begin{figure}[htb]
\center{
\leavevmode
\epsfbox{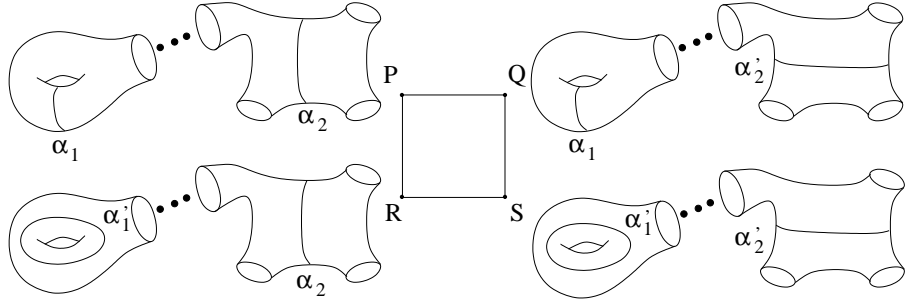}   
\caption{Squares in the Pants Complex.}
\label{squpic}
}
\end{figure}

\begin{figure}[htb]
\center{
\leavevmode
\epsfbox{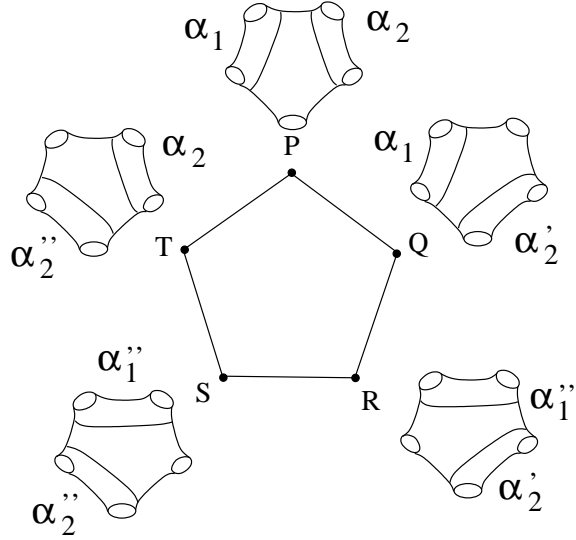}   
\caption{Pentagons in the Pants Complex.}
\label{penpic}
}
\end{figure}

\begin{figure}[htb]
\center{
\leavevmode
\epsfbox{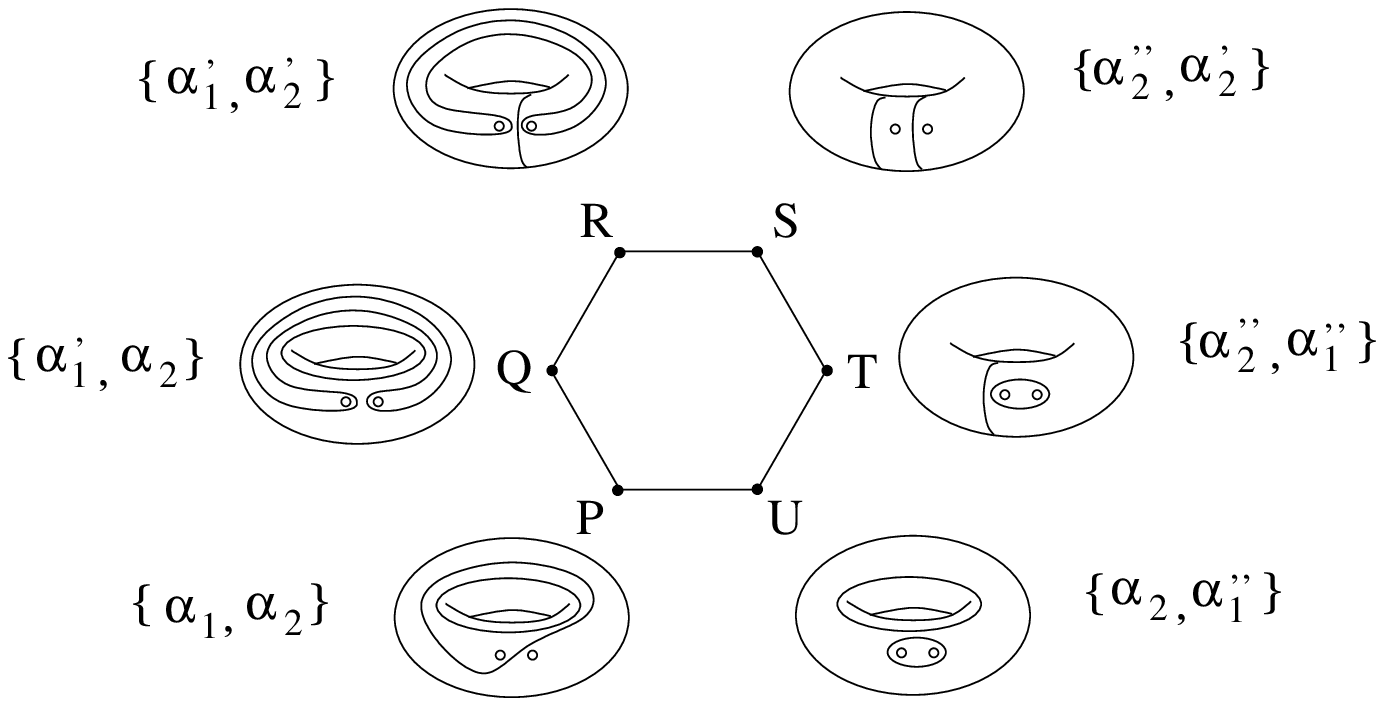}
\caption{Hexagons in the Pants Complex.}
\label{hexpic}
}
\end{figure}

\paragraph{Complex of Pants Decompositions.} The {\em complex of pants decompositions} (or {\em pants complex}) of
$S$, denoted $\cps$, has the pants graph as its $1$-skeleton, and it also has $2$-cells representing specific
relations between elementary moves which are given by topological data on $S$, as depicted in Figures~\ref{tripic}
-~\ref{hexpic}.  In order to prove the main result, we first prove that the pants graph has the same automorphism
group as the pants complex.  Both the pants graph and the pants complex were introduced by Hatcher-Thurston, who used
it to give a presentation for the mapping class group \cite{ht}.  In particular, they showed that the pants complex
is connected and simply connected.

\bigskip

\epsfysize=4 cm

\epsfysize=4 cm

\bigskip

\paragraph{Acknowledgements.}The author would like to thank his advisor, Benson Farb, for suggesting the problem, and
both Benson Farb and Jeff Brock for many helpful conversations.

\section{Sphere With Five Punctures}

We will examine the second isomorphism of the Main Theorem in the case of the 5-times punctured sphere, because it is
the simplest case not classically known.  The techniques involved will be essential in the general case, even though
the method used does not directly generalize.  In the first part of the section, we give topological descriptions of
the pants decompositions corresponding to vertices, edges, triangles, and Farey graphs in the pants complex.  This
will show how to associate a curve on the surface (and hence a vertex of the curve complex) to each Farey graph in
the pants complex.  In the second part, we use this correspondence to define an isomorphism between the automorphism
groups of the pants complex and the curve complex.

\subsection{Pants Graph Objects}

\paragraph{Preliminaries.} First of all, every non-trivial curve on $\szfi$ is {\em 2-separating}, meaning that it
separates $\szfi$ into two components, and that one of the components contains exactly 2 punctures. Hence,
every vertex of $\cc(\szfi)$ can be represented by an arc (really a class of arcs) connecting two
of the punctures. The original curve is recovered by taking the boundary of a small neighborhood
of the arc.

\bigskip

We now build the relationship between Farey graphs in the pants complex and vertices of the curve
complex by first understanding simpler objects in the pants complex.  Since Farey graphs are
preserved by automorphisms of the pants complex, this will give a map from the automorphism group of the
pants complex to the automorphism group of the curve complex. 

\paragraph{Vertices.} Every pants decomposition consists of exactly two curves.  Hence, every vertex
of $\cpo(\szfi)$ can be represented by two arcs on $\szfi$ connecting  two different pairs of punctures. 

\begin{figure}[htb]
\center{
\leavevmode
\epsfbox{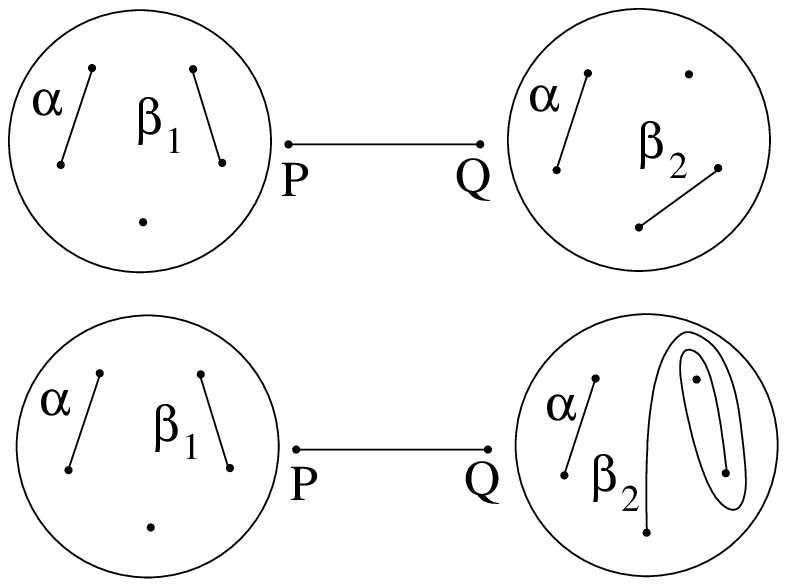}
\caption{Examples of edges in $\cpo(\szfi)$.}
\label{edges}
}
\end{figure}

\paragraph{Edges.} Suppose $P$ and $Q$ are vertices connected by an edge in the pants graph.  Then if
$P$ is represented by the arcs $\{\al,\be_1\}$, $Q$ must be represented by arcs $\{\al,\be_2\}$ for some
$\be_2$ with the property that the curves represented by $\be_1$ and $\be_2$ have intersection
number $2$.  Since a shared endpoint between two arcs corresponds to two points of intersection between
the curves they represent, and a crossing between two arcs corresponds to four points of intersection
between the curves they represent, it follows that there are representatives of $\be_1$ and $\be_2$
which share one puncture as an endpoint and otherwise do not intersect.  The only choice is how many
times $\be_2$ ``winds around'' $\be_1$ (see Figure~\ref{edges}). 

\paragraph{Triangles.} Suppose $P$, $Q$, and $R$ are the vertices of a triangle (a complete graph on three vertices)
in the pants graph. Since the pants decompositions corresponding to $P$ and $Q$ differ by a elementary move, they must
share a common curve.  If $P$ and $Q$ are represented by pairs of arcs $\{\al,\be_1\}$ and
$\{\al,\be_2\}$, then a pair of arcs representing $R$ must have exactly one arc in common with each of
these. This pair of arcs corresponding to $R$ must in fact contain $\al$ (otherwise, it would have to contain
both $\be_1$ and $\be_2$, which can't happen since $\be_1$ and $\be_2$ have nontrivial intersection). Hence $R$ is
represented by $\{\al,\be_3\}$, for some $\be_3$.  In order for $\{\al,\be_3\}$ to differ from $\{\al,\be_1\}$ and
$\{\al,\be_2\}$ by elementary moves, $\be_3$ must have minimal intersection with both $\be_1$ and $\be_2$. This
amounts to finding an arc which shares one endpoint with $\be_1$ and
$\be_2$, and which has no other intersections with $\be_1$, $\be_2$, and $\al$. There
are exactly two such choices. Note that the arcs $\be_1$, $\be_2$, and $\be_3$ (together
with the punctures) form a topological triangle on $\szfi$ which is disjoint from $\al$.

\begin{figure}[htb]
\center{
\leavevmode
\epsfbox{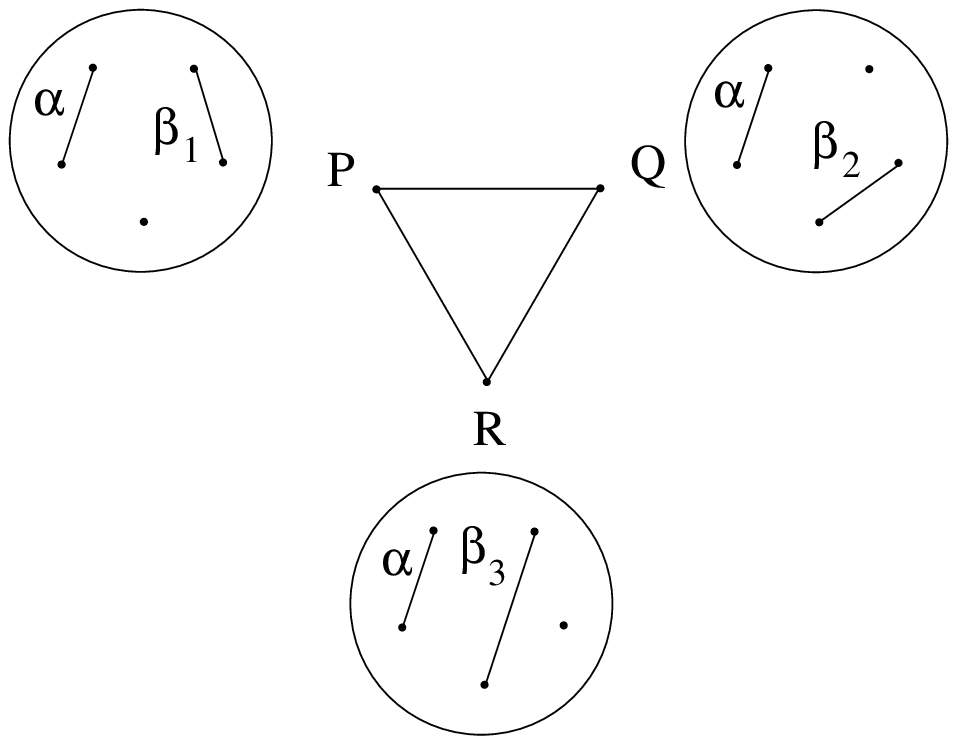}
\caption{A triangle in $\cpo(\szfi)$.}
}
\end{figure}

\paragraph{Farey Graphs.} Any pair of triangles in a single Farey graph can be joined by a path of
triangles in the Farey graph such that consecutive triangles in the path share an edge. Since the pants
decompositions corresponding to any triangle are of the form $\{\al,\be_1\}$, $\{\al,\be_2\}$, and
$\{\al,\be_3\}$, it follows that there is a unique {\em fixed curve} ($\al$) which is in each pants
decomposition corresponding to vertices of the Farey graph in the pants graph. Therefore, there is a
vertex of the curve complex (corresponding to the fixed curve) associated to each Farey graph in
the pants graph.

Likewise, given a vertex $v$ of the curve complex corresponding to a curve $\al$, there is an associated Farey graph
in $\cpo (\szfi)$:  Consider the set of all pants decompositions of $\szfi$ containing the curve $\al$.  This set is
in one-to-one correspondence with the pants decompositions of the $\szfo$ which is a component of the complement of
$\al$ in $\szfi$ and therefore corresponds to a copy of $\cpo(\szfo)$, a Farey graph, denoted $F_v$.

\subsection{The Isomorphism}

\paragraph{Definition of the isomorphism.} In order to prove that $\aut\cpo(\szfi) \cong \mcg(\szfi)$, we
will find an isomorphism $\phi:  \aut\cpo(\szfi) \lra \aut\cc(\szfi)$ and then use the fact that
$\aut\cc(\szfi)  \cong \mcg(\szfi)$.

Given an automorphism $A$ of $\cpo(\szfi)$, we will first define $\phi(A)$, a map from the set of vertices
of $\cc(\szfi)$ to itself. We will then show that $\phi(A)$ is in fact an automorphism of $\cc(\szfi)$, and
that the map $\phi : \aut\cpo(\szfi) \lra \aut\cc(\szfi)$ is an isomorphism.

We define $\phi(A)$ by saying what it does to any particular vertex:  Let $v$ be a vertex of $\cc(\szfi)$,
and $F_v$ the associated Farey graph in $\cpo(\szfi)$. Then $\phi(A)(v)$ is defined to be the vertex of the
curve complex associated to the Farey graph $A(F_v)$, i.e. $A(F_v) = F_{\phi(A)(v)}$. 

\paragraph{$\phi(A)$ is an automorphism of $\cc(\szfi)$.} We will show that $\phi(A)$ takes vertices connected by an
edge to vertices connected by an edge in the curve complex.  It then follows from the definition of the curve complex
that $\phi(A)$ takes the vertices of higher-dimensional simplices to vertices of higher-dimensional simplices (since
every set of $k+1$ mutually connected vertices is the $0$-skeleton of a $k$-simplex), and hence extends to an
automorphism of $\cc(\szfi)$.

Suppose $v$ and $w$ are vertices in the curve complex, corresponding to curves $\al$ and $\be$, and let
$F_v$ and $F_w$ be the Farey graphs associated to $v$ and $w$.  We will show that $v$ and $w$ are connected
by an edge in the curve complex if and only if $F_v$ and $F_v$ have a (unique) common vertex in the pants
complex.  Since the automorphism $A$ preserves the number of intersections between Farey graphs in the pants
complex, this will mean that the map $\phi(A)$ preserves the property of two vertices being connected by an
edge in the curve complex; in other words, $\phi(A)$ extends to an automorphism of the $1$-skeleton of
$\cc(\szfi)$, and hence to an automorphism of $\cc(\szfi)$.

There is a one-to-one correspondence between vertices of $F_v$ and pants decompositions of the form
$\{\al,\star\}$, and a one-to-one correspondence between vertices of $F_w$ and pants decompositions of the
form $\{\be,\star\}$.  Therefore, there is a vertex in common to $F_v$ and $F_w$ if and only if
$\{\al,\be\}$ is a pants decomposition.  But $\{\al,\be\}$ is a pants decomposition exactly when $\al$ and
$\be$ have trivial intersection, which is equivalent to $v$ and $w$ being connected by an edge. 

\paragraph{$\phi$ is an isomorphism.} That $\phi$ is a homomorphism from $\aut\cpo(\szfi)$ to
$\aut\cc(\szfi)$ is clear; it remains to show injectivity and surjectivity.

Injectivity: Suppose that $\phi(A)$ is the identity in $\aut\cc(\szfi)$.  We will show that $A$ is the
identity in $\aut\cpo(\szfi)$.  Since $\phi(A)$ is the identity, $A$ maps every Farey graph in $\cpo(\szfi)$
to itself (given a Farey graph $F_v$ in $\cpo(\szfi)$ corresponding to a vertex $v$ in $\cc(\szfi)$, $A(F_v) 
= F_{\phi(A)v}$ by definition, but $\phi(A)v = v$, so $F_{\phi(A)v} = F_v$). Also, any vertex $P$ is the
unique intersection point of two Farey graphs (if $P$ is $\{\al,\be\}$, where $\al$ and $\be$
correspond to the vertices $v$ and $w$ of the curve complex, then $P = F_v \cap F_w$); hence $A(P)$ must be
equal to $P$. 

Surjectivity: The map $\phi$ has a natural right inverse $\psi: \aut\cc(\szfi) \lra \aut\cpo(\szfi)$, which is
provided by the isomorphism $\eta: \aut\cc(\szfi) \lra \mcg(\szfi)$.  Note that the isomorphism $\eta$ is natural in
the sense that if $\mathcal{A}$ is in $\aut\cc(\szfi)$ and $v$ (corresponding to the curve $\al$) is a vertex of
$\cc(\szfi)$, then $\eta(\mathcal{A})(\al)$ corresponds to the vertex $\mathcal{A}(v)$.

Let $\mathcal{A}$ be an automorphism of the curve complex, and let $P$ be a vertex of the pants graph
which corresponds to the pants decomposition $\{\al,\be\}$.  Then let $\psi(\mathcal{A})(P)$ be the vertex
of the pants graph corresponding to the pants decomposition
$\{\eta(\mathcal{A})(\al),\eta(\mathcal{A})(\be)\}$.  The map $\psi(\mathcal{A})$ is an automorphism of
$\cpo(\szfi)$ because of the natural map $\mcg(\szfi) \lra \aut\cpo(\szfi)$.  It remains to show that
$\phi \circ \psi (\mathcal{A})  = \mathcal{A}$. 

Let $v$ be a vertex of the curve complex.  It suffices to show that if $\mathcal{A}(v) = w$, then $\phi
\circ \psi(\mathcal{A})(v) = w$, but this is by definition the same as $\psi (\mathcal{A}) (F_v) = F_w$. 
The vertices of $F_v$ correspond to all pants decompositions of the form $\{\al, \star\}$, where
$\al$ is the curve corresponding to $v$.  Therefore, by definition, the vertices of $\psi (\mathcal{A})
(F_v)$ correspond to pants decompositions of the form $\{\be, \star\}$, where $\be$ corresponds to
$w$---but this is exactly $F_w$.

We have proven that $\aut\cpo(\szfi) \cong \aut\cc(\szfi)$, so by the fact that $\aut\cc(\szfi) \cong
\mcg(\szfi)$, we are done.

\section{Proof of $\aut \cpos \simeq \aut \cps$}
\label{graph-complex}

In order to prove the first isomorphism of the Main Theorem, we give combinatorial definitions of several kinds of
loops: triangles, alternating squares, alternating pentagons, and almost-alternating hexagons.  These
combinatorial definitions will be stated completely in terms of properties of the pants graph, and so these objects
will be respected by automorphisms of the pants graph.  We prove that our definitions for these loops correspond
exactly to the definitions of the 2-cells in the pants complex given by Figures~\ref{tripic} -~\ref{hexpic}.  This is
enough to show that the automorphism groups of the pants graph and pants complex are canonically isomorphic.  We work
with a fixed surface $S$, which has $n$ curves in each of its pants decompositions (if $S$ is genus $g$ with $r$
punctures, then $n=3g-3+r$).

We first build a correspondence between vertices of the curve complex and \emph{marked} Farey graphs in
the pants graph, which will be used in the sequel to define the map from $\aut\cpos$ to $\aut\cc(S)$.  The definitions
of the various loops will also be used to show that the map between $\aut\cpos$ and $\aut\cc(S)$ is well-defined

A \emph{loop} in the pants graph is a union of a finite collection of edges and the vertices with the property that
each vertex in the collection is the endpoint of exactly two edges in the collection.  We can speak of {\em
consecutive points} in a loop as being a sequence of points such that adjacent points in the sequence share an
edge.  Triangles, squares, pentagons, and hexagons are defined to be loops with the appropriate number of vertices.

Below, a {\em pants graph definition} of an object is a characterization which uses only properties of the pants
graph, and makes no mention of the surface $S$.

\paragraph{Triangles.} {\em We give a pants graph definition of triangles and show that this is equivalent to
the definition of triangular 2-cells in the pants complex.}

Suppose $P$, $Q$, and $R$ are the vertices of a triangle in the pants graph. Since the pants decompositions
corresponding to $P$ and $Q$ differ by an elementary move, they must differ by exactly one curve.  If $P$ and $Q$
correspond to $\{\al_1,\dots,\al_n\}$ and $\{\al_1',\al_2,\dots,\al_n\}$, then a pants decomposition corresponding to
$R$ must have exactly $n-1$ curves in common with each of these.  This pants decomposition associated to $R$ must in
fact contain $\al_2,\dots,\al_n$ (otherwise, it would have to contain $\al_1$ and $\al_1'$, which can't happen since
$\al_1$ and $\al_1'$ intersect).  Hence $R$ corresponds to $\{\al_1'',\al_2,\dots,\al_n\}$, for some $\al_1''$.  In
order for $\{\al_1'',\al_2,\dots,\al_n\}$, to differ from $\{\al_1,\dots,\al_n\}$ and $\{\al_1',\al_2,\dots,\al_n\}$
by elementary moves, $\al_1''$ must have minimal intersection with both $\al_1$ and $\al_1'$.

There are two topologically distinct possibilities for triangles in the pants graph, depending on
whether $P$ and $Q$ differ by an elementary move on a $\szfo$, or an elementary move on a $\soo$.  In the first
case, $\al_1$, $\al_1'$, and $\al_1''$ all lie on the $\szfo$, and in the second case they lie on the
$\soo$.  So understanding triangles on a general surface boils down to understanding triangles on these
subsurfaces. Since the pants complex for both the sphere with four punctures and the torus with one
puncture is a Farey graph, and $\al_1$ and $\al_1'$ are connected by an edge in these complexes, there
are exactly two choices for $\al_1''$ in either case.

Note that the two types of triangles described here correspond exactly to the two types of triangular
2-cells in the pants complex given by Figure~\ref{tripic}.

\paragraph{Marked Farey Graphs.} {\em We give a correspondence between marked Farey graphs in the pants graph and
vertices in the curve complex (to be used in the definition of the isomorphism between $\aut \cps$ and $\aut \cc(S)$
in the next section).}

By the chain-connectedness of Farey graphs (any two triangles can be connected by a sequence of triangles so that
consecutive triangles in the sequence share an edge), and since the pants decompositions corresponding to any
triangle are of the form $\{\al_1^1,\al_2,\dots,\al_n\}$, $\{\al_1^2,\al_2,\dots,\al_n\}$, and
$\{\al_1^3,\al_2,\dots,\al_n\}$, it follows that there are $n-1$ {\em fixed curves} ($\al_2,\dots,\al_n$), and one
{\em moving curve} (the $\al_1^i$'s) in the pants decompositions corresponding to vertices of the Farey graph. By
marking a point on a Farey graph, we are essentially marking one of the $\al_1^i$'s.  Hence, there is a unique vertex
of the curve complex associated to a marked Farey graph in the pants graph.

Likewise, given a vertex $v$ of the curve complex corresponding to a curve $\al_1$, there is an
associated marked Farey graph in the pants graph---but it is not unique.  To get a marked Farey graph
associated to $v$, construct a pants decomposition containing $\al_1$, say $\{\al_1,\al_2,\dots,\al_n\}$,
and consider the set of all pants decompositions containing $\al_2,\dots,\al_n$.  Since the complement of
$\al_2,\dots,\al_n$ in the surface is a number of pants and either a $\szfo$ or
$\soo$, the set of pants decompositions of the form $\{\star,\al_2,\dots,\al_n\}$ corresponds to a Farey
graph in the pants graph (because both $\cpo(\szfo)$ and $\cpo(\soo)$ are Farey graphs). With the vertex
corresponding to $\{\al_1,\al_2,\dots,\al_n\}$ marked, we have a marked Farey graph corresponding to
$\al_1$.  A marked Farey graph representing a vertex $v$ of the curve complex is denoted $(F_v,P)$, where
$F_v$ is a Farey graph, and $P$ is the marked vertex. 

\paragraph{Definition of Alternating Sequences.} {\em We give a pants graph definition of alternating sequences,
which is a key concept used in the next section to show that the isomorphism between $\aut \cps$ and $\aut \cc(S)$ is
well-defined.}

In our discussions of squares, pentagons, and hexagons, we will use two equivalent notions of an \emph{alternating}
sequence of points in a loop in the pants graph---one geometric and one combinatorial.

The sequence of consecutive points $P_1P_2\dots P_n$ in a loop $\lp$ is alternating if for any three
consecutive points $P_iP_{i+1}P_{i+2}$ in the sequence

\begin{itemize}

\item there is no Farey graph in the pants graph containing $P_iP_{i+1}P_{i+2}$ (or the pair of edges
connecting these points) (geometric definition)

or, equivalently

\item the pants decompositions corresponding to $P_iP_{i+1}P_{i+2}$ have no set of $n-1$ curves in
common. (combinatorial definition)

\end{itemize}

That these definitions are equivalent follows from the characterization of Farey graphs above.  A useful
working definition of an alternating sequence of points $PQR$ is that if the elementary move corresponding
to $PQ$ is $\star \ra \al$, then the elementary move corresponding to $QR$ is not of the form $\al \ra
\star$. Any loop in the pants graph with the property that any three consecutive vertices make up an
alternating sequence is called an \emph{alternating loop}. 

\paragraph{Alternating Squares.} {\em We give a pants graph definition of alternating squares and show that this is
equivalent to the definition of square 2-cells in the pants complex.}

An \emph{alternating square} in the pants graph is a square which is alternating.  By the combinatorial definition of
alternating, it follows that pants decompositions corresponding to opposite corners of the square differ by two
curves.  So an alternating square is a ``true square'' in the sense that it has no extra connections.  The converse
statement, that all true squares are alternating squares, will follow from the argument below.

Suppose $P$, $Q$, $R$ and $S$ are the (ordered) vertices of an alternating square in the pants graph,
where $P$ and $Q$ correspond to the pants decompositions $\{\al_1,\al_2, \dots, \al_n\}$ and
$\{\al_1,\al_2',\dots, \al_n\}$ for some $\al_2'$ which has minimal intersection with $\al_2$. In order for
$R$ to be connected to $Q$ and to not be in the same Farey graph as $P$ and $Q$, its pants decomposition
must contain $\al_2'$, and so must be of the form $\{\al_1',\al_2',\dots, \al_n\}$ for some $\al_1'$ which
has minimal intersection with $\al_1$ (since the elementary move corresponding to $PQ$ is $\al_2 \ra \al_2'$,
the elementary move corresponding to $QR$ cannot be of the form $\al_2' \ra \star$). Now $S$ must be connected
to both $P$ and $R$. Therefore the pants decomposition associated to $S$ must contain $\al_3,\dots,\al_n$
(because it can only differ from each of the others by one curve); also, it must not contain $\al_1$
(otherwise $S$, $P$, and $Q$ would lie in a Farey graph); it must also contain $\al_2'$ (so as to differ
from the pants decomposition associated to $R$ by one curve.  Therefore, the pants decomposition
corresponding to $S$ must be $\{\al_1',\al_2,\al_3,\dots, \al_n\}$

We can be more specific about what these pants decompositions look like on the surface.  In particular,
we will see that the pair $\al_1, \al_1'$ is on a disjoint subsurface from the pair $\al_2,\al_2'$ in
the complement of $\al_3,\dots,\al_n$.  This will follow from a more general
principle, stated and proven below, that if the edges $PQ$ and $QR$ in the pants graph form two sides of
an alternating square, where $PQ$ and $QR$ correspond to moves $\al_1 \ra \al_1'$ and $\al_2 \ra \al_2'$
for some curves $\al_1$, $\al_1'$, $\al_2$ and $\al_2'$, then $\al_1$ and $\al_2$ lie on disjoint
subsurfaces.

Note that the alternating squares described here correspond exactly to the square 2-cells in the
pants complex given by Figure~\ref{squpic}.

\paragraph{Half-Squares.} {\em We give a pants graph definition of half-squares and use it to understand
alternating squares and 3-curve small loops (see below).}

A \emph{half-square} in the pants graph is a set of three vertices, say $P$, $Q$, and $R$, along with two edges $PQ$
and $QR$ with the property that there is a fourth vertex in the pants graph, say $S$, so that $P$, $Q$, $R$, and $S$
are the vertices of an alternating square in the pants graph. If $PQ$ and $QR$ correspond to elementary moves $\al_2
\ra \al_2'$ and $\al_1 \ra \al_1'$ for some curves $\al_1$, $\al_1'$, $\al_2$ and $\al_2'$, then $P$, $Q$, and $R$
are the vertices of a half-square if and only if $\al_1'$ and $\al_2$ have trivial intersection; for if this is the
case, then edges corresponding to moves $\al_2' \ra \al_2$ and $\al_1 \ra \al_1'$ will ``complete the square''. This
clearly happens when $\al_1$ and $\al_2$ lie on disjoint subsurfaces, but the claim is that this condition is
necessary.

Assume that $P$, $Q$, $R$, are the vertices of a half-square, and that $\al_1$ and $\al_2$ lie on the same subsurface,
either a $\szfi$ or a $\sot$.  There are four topological possibilities (up to the action of the extended mapping
class group for either of these surfaces) for $\al_1$, $\al_2$, and $\al_2'$---on the $\szfi$ there is only one
possibility, and on the $\sot$ there are three cases (since any pants decomposition of $\sot$ contains at most one
separating curve, and since an elementary move always changes a separating curve to a nonseparating one, there can be
at most one separating curve in a half-square, so the cases correspond to $0$, $1$, and $2$ appearances of a
separating curve): \emph{(i)} $\al_1$, $\al_2$, and $\al_2'$ are all nonseparating, \emph{(ii)} $\al_2'$ and $\al_1$
are nonseparating and $\al_2$ is separating, and \emph{(iii)} $\al_2$ and $\al_2'$ are nonseparating and $\al_1$ is
separating.  It is clear that in each of the cases, there is no curve $\al_1'$ which intersects $\al_1$ minimally and
is disjoint from $\al_2$ and $\al_2'$.  This is a contradiction, so the result follows.

\paragraph{Small Loops.} {\em We give a pants graph definition of small loops and give a characterization of the
corresponding pants decompositions; we use this to understand alternating pentagons and almost-alternating hexagons
below.}

A {\em small loop} in the pants graph is a loop with no more than six edges.  It turns out that for any small loop
$\lp$ in a pants graph, exactly one of the following is true:

\emph{(a)} the pants decompositions corresponding to the vertices of the loop all have the same $n-2$
curves in common, i.e. they are of the form $\{\star,\star,\al_3,\dots,\al_n\}$.  In this case,
$\lp$ is called a {\em 2-curve small loop}.  Note that an alternating square is a 2-curve small loop. 
Alternating pentagons and almost-alternating hexagons will also be 2-curve small loops. 

\emph{(b)} $\lp$ is a hexagon, and the pants decompositions corresponding to the vertices of the loop all
have the same $n-3$ curves in common, i.e. they are of the form
$\{\star,\star,\star,\al_4,\dots,\al_n\}$, and any pair of edges in the loop which share an endpoint are two
sides of an alternating square in the pants graph.  In this case, $\lp$ is
called a {\em 3-curve small loop}. 

The thrust of this statement is that if there is a small loop in the pants graph, then either all but
two curves in the corresponding pants decompositions stay fixed, or something else very specific happens.

We will show that if there is no set of $n-2$ curves that appear in each pants decomposition
corresponding to the vertices of $\lp$ (i.e. if $\lp$ is not a 2-curve small loop), then the loop $\lp$
is a hexagon which has the property that consecutive edges are the edges of an alternating square.

Suppose there is a vertex of $\lp$ corresponding to the pants decomposition $\{\al_1,\dots,\al_n\}$. Recall
that (directed) edges of $\lp$ represent elementary moves, which can be denoted, for example, $\al_1 \ra
\al_1'$.  The directions of the edges should match up with a given orientation of $\lp$.  Assume without
loss of generality that there is at least one vertex of $\lp$ that does not contain $\al_1$, at least one
that does not contain $\al_2$, and at least one that does not contain $\al_3$ (i.e. assume that there is not
a set of $n-2$ curves which stay fixed).

In this case, there must be three edges in $\lp$ corresponding to the moves $\al_1 \ra \star$, $\al_2
\ra \star$, and $\al_3 \ra \star$ (after having chosen a preferred direction around the loop).  In order
for $\lp$ to be a loop, it must also have edges corresponding to moves $\star \ra \al_1$, $\star \ra
\al_2$, and $\star \ra \al_3$.  These two sets of moves must be distinct---none of the $\star$'s can be
$\al_1$, $\al_2$, or $\al_3$, since these curves have trivial intersection pairwise.  Hence, the
moves of $\lp$ can be written as $\al_1 \ra \al_1'$, $\al_2 \ra \al_2'$, $\al_3 \ra \al_3'$, $\star \ra
\al_1$, $\star \ra \al_2$, and $\star \ra \al_3$, (the $\al_i'$ are not assumed to be distinct) and
$\lp$ must have at least six edges. 

Assume that $\lp$ has no more edges than the ones already described, i.e. that $\lp$ has no more than six
edges.  Then the three moves $\star \ra \al_1$, $\star \ra \al_2$, and $\star \ra \al_3$ must be in fact
$\al_1' \ra \al_1$, $\al_2' \ra \al_2$, $\al_3' \ra \al_3$.  By investigating the possible orders of these
six edges, we will determine that exactly one of the following is true.

\emph{(i)} There must be an edge outside of $\lp$ connecting two vertices of $\lp$ (i.e. $\lp$ is not a true
loop).

\emph{(ii)} The edges of $\lp$ are, in order, $\al_1 \ra \al_1'$, $\al_2 \ra \al_2'$, $\al_3 \ra
\al_3'$, $\al_1' \ra \al_1$, $\al_2' \ra \al_2$, $\al_3' \ra \al_3$.  In this case, $\lp$ is a squared
hexagon (see the characterization of small loops).

Assume that \emph{(ii)} does not happen.  We will show that \emph{(i)} happens, i.e. $\lp$ is not a true
loop. Since we are assuming it is not the case that the edges corresponding to $\al_i \ra \al_i'$ and
$\al_i' \ra \al_i$ are on opposite sides of the hexagon, then it must be the case that in one direction
around $\lp$, there are three consecutive edges $e_1$, $e_2$, and $e_3$ corresponding to moves $\al_i \ra
\al_i'$, $\al_j \ra \al_j'$, and $\al_i' \ra \al_i$.  But this means that $\al_j'$ has trivial intersection
with $\al_i$, and so the edges corresponding to the moves $\al_i \ra \al_i'$ and $\al_j \ra \al_j'$ are the
edges of a half-square whose other two sides correspond to moves $\al_i' \ra \al_i$ (this is edge $e_3$) and
$\al_j' \ra \al_j$ (call this edge $e_4$).  Since $e_4$ connects endpoints of $e_1$ and $e_3$, it follows
that $\lp$ is not a true loop---it is two squares sharing an edge. 

In the case that \emph{(ii)} happens, we will show that pairs of consecutive edges are the edges of
half-squares in the pants graph.  By the characterization of half squares, this is the same as $\al_i$
and $\al_j'$ having trivial intersection for any $i$ and $j$.  But this is the case, since (for any $i$
and $j$) $\al_i$ and $\al_j'$ appear in a common pants decomposition corresponding to some vertex of the
loop.

\paragraph{Alternating Pentagons.} {\em We give a pants graph definition of alternating pentagons and show that this
is equivalent to the definition of pentagonal 2-cells in the pants complex.}

An \emph{alternating pentagon} in the pants graph is a pentagon which is alternating, i.e. it has the property that
no pair of consecutive sides lie in a single Farey graph. Again, by the combinatorial definition of alternating, this
means that the pentagon does not have three consecutive vertices whose pants decompositions share $n-1$ curves in
common.

Let $P$, $Q$, $R$, $S$, and $T$ be the (ordered) vertices of an alternating pentagon in the pants graph. 
By the characterization of small loops (a pentagon is a small loop), the pants decompositions corresponding
to these vertices all have $n-2$ curves in common, say $\al_3,\dots,\al_n$.  If the pants
decompositions corresponding to $P$ and $Q$ are $\{\al_1,\al_2,\dots,\al_n\}$ and $\{\al_1,\al_2',\dots,\al_n\}$
for some $\al_1$, $\al_2$, and $\al_2'$, then the vertex corresponding to $R$ must be
$\{\al_1',\al_2',\dots,\al_n\}$ for some $\al_1'$ (otherwise $P$, $Q$, and $R$ all have $n-1$ curves
in common, contradicting the fact that the pentagon is alternating).  Since $S$ and $T$ are connected to $R$
and $P$, respectively, then their associated pants decompositions must be of the form
$\{\al_1',\al_2'',\dots,\al_n\}$ and $\{\al_1'',\al_2,\dots,\al_n\}$ for some $\al_1''$ and $\al_2''$ (the
pants decompositions corresponding to $S$ cannot contain $\al_2'$, because then the pants decompositions
corresponding to $Q$, $R$, and $S$ would have $n-1$ curves in common;  likewise for $T$).  Now, since the
vertices $S$ and $T$ are connected by an edge, their associated pants decompositions,
$\{\al_1',\al_2'',\dots,\al_n\}$ and $\{\al_1'',\al_2,\dots,\al_n\}$, must have $n-1$ curves in
common. Also, since the edge $RS$ corresponds to the move $\al_2' \ra \al_2''$, the edge $ST$ must
correspond to a move $\al_1' \ra \star$ (it can't be $\al_2'' \ra \star$ because that would violate the
alternating property).  Likewise, since the move corresponding to $PT$ is $\al_1 \ra \al_1''$, the move
corresponding to $TS$ must be of the form $\al_2 \ra \star$.  Putting the characterizations of $ST$ and $TS$
together, the move corresponding to $ST$ must be $\al_1' \ra \al_2$.  But it then follows that $\al_2''$ and
$\al_1''$ are the same curve (by comparing the lists of curves in the pants decompositions
corresponding to $S$ and $T$). 

In summary, the pants decompositions corresponding to $P$, $Q$, $R$, $S$, and $T$ are
$\{\al_1,\al_2,\dots,\al_n\}$, $\{\al_1,\al_2',\dots,\al_n\}$, $\{\al_1',\al_2',\dots,\al_n\}$,
$\{\al_1',\al_2'',\dots,\al_n\}$, $\{\al_2,\al_2'',\dots,\al_n\}$ and the elementary moves corresponding to
$PQ$, $QR$, $RS$, $ST$, and $TP$ are $\al_2 \ra \al_2'$, $\al_1 \ra \al_1'$, $\al_2' \ra \al_2''$, $\al_1'
\ra \al_2$, and $\al_2'' \ra \al_1$.  Using the facts that any two curves in the same pants
decomposition are disjoint, and that any curves $\be$ and $\ga$ making up an elementary move $\be \ra \ga$
intersect minimally, we have that in the sequence of curves $\al_2$, $\al_2'$, $\al_2''$, $\al_1$,
$\al_1'$, $\al_2$, curves which are adjacent in the sequence intersect minimally, and curves
which are not adjacent in the sequence are disjoint.

First of all, this implies that $\al_1$ and $\al_2$ do not lie on disjoint subsurfaces (since $\al_1'$
 has nontrivial intersection with both of them).  Therefore, $\al_1$ and $\al_2$ must lie on a $\szfi$ or $\sot$ in
the complement of $\al_3,\dots,\al_n$.  In the first case, the curves $\al_2$, $\al_2'$,
$\al_2''$, $\al_1$, and $\al_1'$ must be as in the definition of pentagonal 2-cells in the pants complex (Figure 5),
up to the action of the mapping class group.  We now show that the second case cannot happen, i.e.  that there is no
such sequence of five curves on a torus with two punctures.

Assume that on $\sot$ there is a sequence of curves $\al$, $\be$, $\gamma$, $\delta$, $\epsilon$,
$\al$ with the property that consecutive curves intersect minimally and curves that are
not adjacent in the sequence do not intersect.  In such a sequence, there can be at most one curve
which is separating on $\sot$.  This is because two separating curves on $\sot$ intersect at
least four times (since any separating curve can be obtained by taking the boundary of a
small neighborhood of an arc connecting the two punctures), but any two curves in the sequence
must intersect no more than two times since they are either disjoint or minimally intersecting.  We will
consider two cases, corresponding to whether there is a separating curve in the sequence (Case 1)
or there is no separating curve in the sequence (Case 2). Both cases will result in a
contradiction. We call a nonseparating curve on $\sot$ of $(p,q)$-type on if it is of $(p,q)$ class on the
torus obtained by ``forgetting'' the two punctures. 

\case{1} Suppose there is a separating curve in the sequence, say $\al$.  It follows that the
other curves in the sequence are nonseparating and that $\al$ separates
$\sot$ into a pair of pants and a punctured torus. Since $\gamma$ and $\delta$ both have trivial
intersection with $\al$ and have minimal intersection with each other in the complement of a
$\al$, they must lie on the punctured torus and intersect once.  Say that $\gamma$ is
of $(1,0)$-type, and $\delta$ is of $(0,1)$-type, and so a curve is of type $(p,q)$ if it
intersects $\delta$ $p$ times and $\gamma$ $q$ times.  Since $\be$ and $\epsilon$ are both nonseparating
curves on $\sot$, and since $\be$ has trivial intersection with $\delta$, it must be of
$(0,1)$-type; likewise $\epsilon$ must be of $(1,0)$-type.  This implies that $\be$ and $\epsilon$ have
nontrivial intersection, which is a contradiction.

\begin{figure}[htb]
\center{
\leavevmode
\epsfbox{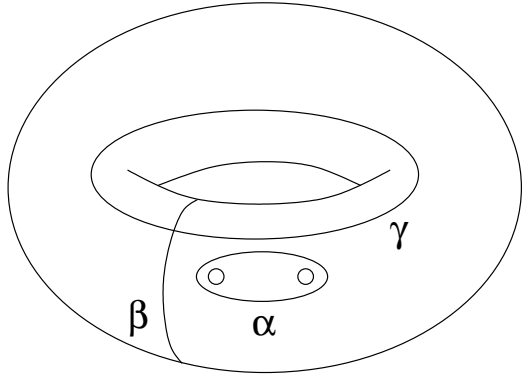}
\caption{Case 1.  The configuration for $\al$ separating.}
}
\end{figure}

\case{2} Now suppose that all the curves in the sequence are nonseparating.  In this case, the
complement of any of the curves is $\szfo$, so a minimal intersection number
between two curves will always be $2$.  Since all pants decompositions of $\sot$ consisting of
two nonseparating curves are equivalent up to the extended mapping class group, we can assume
without loss of generality that $\al$ and $\gamma$ are the two curves shown in the figure. 

\begin{figure}[htb]
\center{
\leavevmode
\epsfbox{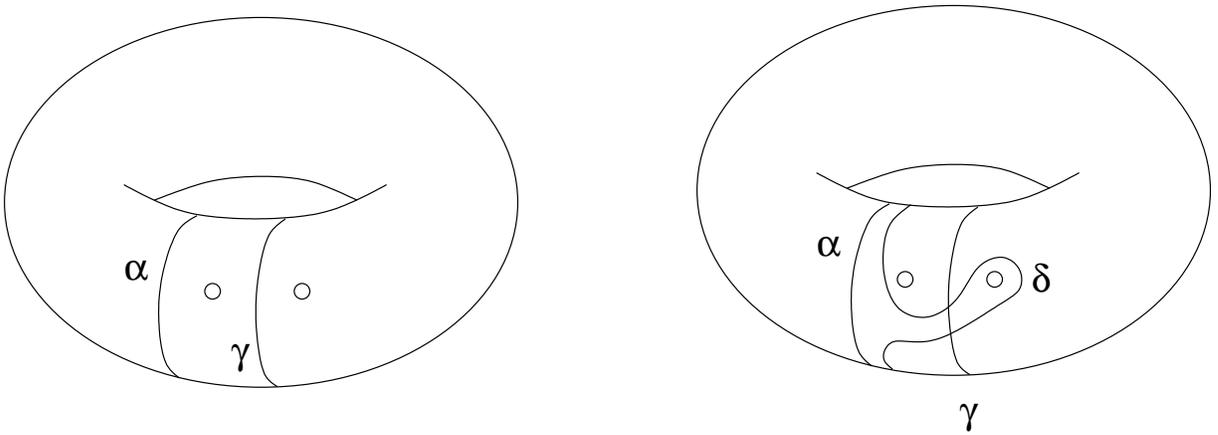}
\caption{Case 2.  (i) The configuration for $\al$ nonseparating.  (ii) The unique choice of $\delta$.} }
\end{figure}

Now since $\delta$ differs from $\gamma$ by an elementary move it cannot have the same association as $\ga$ on the
$\szfo$ which is the the complement of $\al$.  Since it is nonseparating, it cannot associate the
two punctures.  That leaves one possible association for $\delta$, so all choices for $\delta$ are equivalent up to
Dehn twists around $\gamma$ and a reflection (both of which fix $\al$ and $\ga$), by the associativity move rule (see
Section~\ref{complexes}).

Now we will show that there is no curve $\ep$ on $\sot$ which differs by an elementary move from both
$\al$ and $\delta$ and has trivial intersection with $\ga$.  The curves $\ga$ and $\delta$ divide $\sot$
into regions, two of which are punctured disks (the punctures being the two punctures of $\sot$).  The
claim is that any curve $\ep$ which differs from $\al$ by an elementary move on $\sot$ and is disjoint from
$\ga$ must pass in and out of each of these punctured disks through the curve $\delta$, and therefore
has intersection number at least $4$ with $\delta$, and thus does not differ from $\delta$ by an elementary
move on $\sot$.  The proof of the claim is that if we cut $\sot$ along $\ga$ and remove one of the
punctured disks, what is left over is a pair of pants whose three boundaries are the not-removed
puncture from $\sot$, $\ga$, and a curve homotopic to $\al$; any curve on a pair of pants must be
homotopic to one of the boundary components and hence does not differ from $\al$ by an elementary move.  We
will refer to this argument in later sections as the {\it Punctured Disks Argument}.

Since there is no such curve $\ep$, there is no appropriate configuration for $\al$, $\be$, $\gamma$,
$\delta$, and $\epsilon$ on $\sot$, and hence no alternating pentagons on $\sot$. 

Note that the alternating pentagons described here correspond exactly to the pentagonal 2-cells in the
pants complex given by Figure~\ref{penpic}.

\paragraph{Quadrilateral Triples.} {\em We give a pants graph definition of quadrilateral triples which will be used
in the definition of almost-alternating hexagons below.}

Besides the concept of alternating sequences, we need to understand one further type of geometric configuration in
the pants graph before discussing hexagons.  A \emph{quadrilateral triple} of points in the pants graph is a set of
three points in the pants graph which lie in a common quadrilateral in a Farey graph, and which do not lie on a
common triangle.  The unique point in the triple connected to the other two by edges is called the \emph{central
point}, while the other two points are called \emph{outer points}.  In the Farey graphs corresponding to the pants
graph of the punctured torus and the four-times punctured sphere, an example of a quadrilateral triple is the
vertices corresponding to the $(1,0)$ curve, the $(0,1)$ curve, and the $(2,1)$ curve.  Note that the vertex
corresponding to the $(1,1)$ curve is connected to each of these vertices by an edge.  Also, up to the action of the
extended mapping class group on the pants graph for either of these surfaces, all quadrilateral triples are
equivalent.

\paragraph{Almost-Alternating Hexagons.} {\em We give a pants graph definition of almost-alternating hexagons and
show that this is equivalent to the definition of hexagonal 2-cells in the pants complex.}

An \emph{almost-alternating hexagon} is a hexagon in the pants graph which has an alternating sequence
of six vertices, and three vertices which make up a quadrilateral triple. 

Let $P$, $Q$, $R$, $S$, $T$, and $U$ be the (consecutive) vertices of an almost-alternating hexagon,
where $UPQ$ is the quadrilateral triple.  Then the alternating sequence must be $PQRSTU$ (the
quadrilateral triple cannot be a subsequence of the alternating sequence).

Since an almost-alternating hexagon is not a 3-curve small loop, then by the characterization of small
loops in the pants graph, the pants decompositions corresponding to the vertices all have a set of $n-2$
curves in common, say $\al_3,\dots,\al_n$.

If the curves in the pants decomposition corresponding to $P$ are $\{\al_1,\al_2,\dots,\al_n\}$, then
those corresponding to $Q$ and $U$ must be $\{\al_1',\al_2,\dots,\al_n\}$ and
$\{\al_1'',\al_2,\dots,\al_n\}$ for some $\al_1'$ and $\al_1''$ (since they lie in a common Farey graph). 
Then, since the sequence $PQR$ is alternating, the pants decomposition corresponding to $R$ must be
$\{\al_1',\al_2',\dots,\al_n\}$ for some $\al_2'$.  Likewise, the pants decomposition for $T$ must be
$\{\al_1'',\al_2'',\dots,\al_n\}$ for some $\al_2''$.  Finally, since the sequence $QRS$ is alternating, the
pants decomposition corresponding to $S$ must contain $\al_2'$, and since $STU$ is alternating, the pants
decomposition corresponding to $S$ must contain $\al_2''$. Therefore, the pants decomposition corresponding to $S$
must be $\{\al_2',\al_2'',\dots,\al_n\}$.

In order to carry out the topological characterization of almost-alternating hexagons, we take the
following steps: \\
1.  $\al_1$ and $\al_2$ do not lie on disjoint subsurfaces \\
2.   $\al_1$ and $\al_2$ do not lie on a $\szfi$ (and hence they lie on a $\sot$) \\
3.  $\al_2$ is nonseparating on the $\sot$ \\
4.  $\al_1$ is nonseparating on the $\sot$ \\
5.  $\al_1'$ (and hence $\al_1''$) is separating on the $\sot$ \\ 
6.  The choices of $\al_1$, $\al_1'$, $\al_1''$, $\al_2$, $\al_2'$, and $\al_2''$ are unique up to the
action of $\mcg(\sot)$

\step{1} The curves $\al_1$ and $\al_2$ cannot lie on disjoint subsurfaces, since there is a chain
of curves connecting them which are disjoint from $\al_3,\dots,\al_n$.  Namely, $\al_1$ has
intersection with $\al_1''$, which has intersection with $\al_2'$, which has intersection $\al_2$ (all of
these curves have trivial intersection with $\al_3,\dots,\al_n$ since they make up pants
decompositions with these curves). 

\step{2} Assume $\al_1$ and $\al_2$ lie on a $\szfi$ in the complement of
$\al_3,\dots,\al_n$. Since $P$, $Q$, and $U$ are the vertices of a quadrilateral in a Farey graph, and
$\al_2$ appears in all three corresponding pants decompositions, the aforementioned Farey graph is the pants
complex of the $\szfo$ which is one of the components of $\szfi - \al_2$.  Then since any two quadrilaterals
in $\cpo(\szfo)$ are equivalent up to the extended mapping class group of $\szfo$, and any mapping class of
$\szfo$ extends to a mapping class of $\szfi$, it follows that up to the extended mapping class group, the
pants decompositions corresponding to $P$, $Q$, and $U$ are as follows:

\begin{figure}[htb]
\center{
\leavevmode
\epsfbox{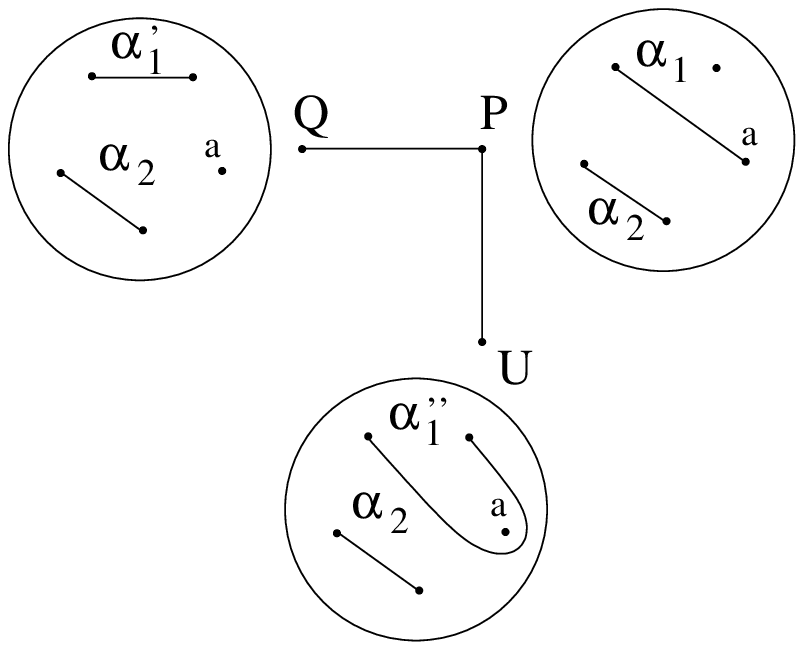}
\caption{Step 2: The configuration for $\al_1$ and $\al_2$ on a $\szfi$.}
}
\end{figure}

Note that we use the arc representations of curves, as discussed in the treatment of $\cpo(\szfi)$
above.  The edges $QR$ and $UT$ (note directions) correspond to the elementary moves $\al_2 \ra \al_2'$ and
$\al_2 \ra \al_2''$.  By the characterization of elementary moves between pants decompositions of $\szfi$, both
$\al_2'$ and $\al_2''$ must be represented by arcs which have an endpoint at the puncture $a$ (see figure).
This implies that $\al_2'$ and $\al_2''$ have nontrivial intersection---a contradiction, since they both
appear in the pants decomposition corresponding to $S$.

Now we know that $\al_1$ and $\al_2$ lie on a $\sot$. 

\step{3} If we assume $\al_2$ is separating on the $\sot$, then it separates $\sot$
into a pair of pants and a punctured torus.  Then $\{\al_1\}$, $\{\al_1'\}$, and $\{\al_1''\}$ are pants
decompositions of the $\soo$, whose corresponding vertices in $\cpo(\soo)$ make a quadrilateral triple. 
Since all quadrilaterals in $\cpo(\soo)$ are equivalent up to $\mcg(\soo)$, and since every mapping class of
$\soo$ extends to a mapping class of $\sot$ which fixes $\al_2$, we assume that on $\soo$, $\al_1''$ is the
$(1,0)$ curve, $\al_1$ is the $(2,1)$ curve, and $\al_1'$ is the $(1,1)$ curve.  These
curves are of $(1,0)$-type, $(2,1)$-type, and $(1,1)$-type on $\sot$ (as above, a nonseparating curve
 in $\sot$ is of $(p,q)$-type if it is the $(p,q)$ curve on the torus obtained by forgetting the
two punctures).  Since $\al_2$ is separating, $\al_2'$ and $\al_2''$ must both be nonseparating (since they
must have intersection number with $\al_2$ no more than $2$). Also, because $\al_2'$ must have trivial
intersection with $\al_1'$ (the two curves form a pants decomposition), it follows that $\al_2'$ must be of
type $(1,1)$.  Likewise, $\al_2''$ must be of type $(1,0)$.  However, since $\al_2'$ and $\al_2''$ must make
up a pants decomposition of $\sot$, they must have trivial intersection; but curves of different type
always have nontrivial intersection.  We have a contradiction, and so $\al_2$ must be nonseparating.


Since $\al_2$ is nonseparating on $\sot$, it follows that the complement of $\al_2$ in
$\sot$ is a $\szfo$, and therefore $\{\al_1'\}$, $\{\al_1\}$, and $\{\al_1''\}$ are pants decompositions of
the $\szfo$.  Note that two of the boundary components in the $\szfo$ correspond to punctures in the $\sot$
($p_1$ and $p_2$), while the other two correspond to the two ``sides'' of $\al_2$
($\al_2^+$ and $\al_2^-$).

\step{4} Assume that $\al_1$ is separating on the $\sot$.  Since all pants decompositions containing a
separating curve are equivalent up to the extended mapping class group, we assume that $\al_1$ and
$\al_2$ are as in the figure.

\epsfysize=4 cm
\begin{figure}[htb]
\center{
\leavevmode
\epsfbox{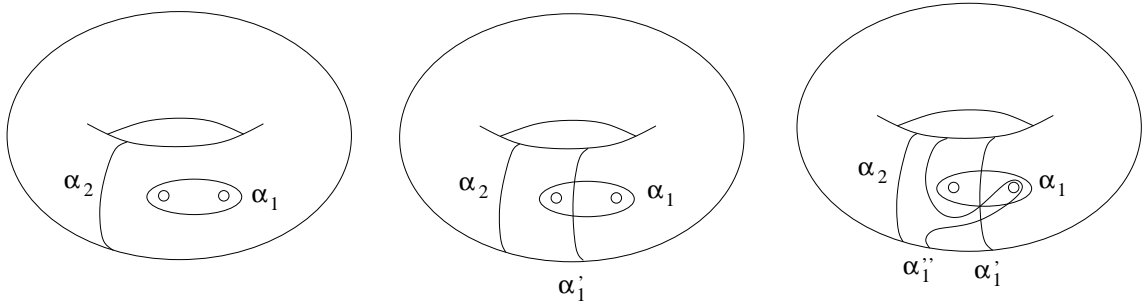}
\caption{Step 4. (i) The configuaration for $\al_1$ separating. (ii) The unique choice for $\al_1'$.
(iii) The unique choice for $\al_1''$.}
}
\end{figure}

Since $\al_1'$ differs from $\al_1$ by an elementary move on the $\szfo$ which is the complement
of $\al_2$, and since $\al_1$ associates the two punctures (it is separating), it follows that $\al_1$ associates a
puncture with one side of $\al_2$.  Up to a mapping class which permutes the punctures and fixes $\al_1$ and $\al_2$,
we can choose which puncture gets associated with which side of $\al_2$.  Therefore, the choice of
$\al_1'$ is unique up to Dehn twists about $\al_1$ and a reflection (both of which fix $\al_1$ and $\al_2$), by the
associativity move rule (see Section~\ref{complexes}).

Since $\al_1''$ must have the same association as $\al_1'$, and it must be part of a quadrilateral triple
with the vertices of $\cpo(\szfo)$ corresponding to $\al_1$ and $\al_1'$, the choice of $\al_1''$ is unique
up to a reflection in $\sot$ which fixes $\al_1$, $\al_2$, and $\al_1'$. 

The curve $\al_2'$ must have trivial intersection with $\al_1'$, and must differ from both
$\al_1''$ and $\al_2$ by elementary moves.  By the Punctured Disks Argument (see {\it Alternating
Pentagons}), $\al_2'$ must have intersection number greater than $2$ with $\al_1''$ and therefore cannot
differ from $\al_1''$ by an elementary move.  This is a contradiction, so $\al_1$ must be nonseparating.




For Steps 5 and 6, we assume that $\al_1$ and $\al_2$ are of $(0,1)$-type.

\step{5} Assume that $\al_1'$ is nonseparating on $\sot$.  Since $\al_1$ is nonseparating, and $\al_1'$ differs from
$\al_1$ by an elementary move on the $\szfo$ which is the complement of $\al_2$ on $\sot$, there
is only one possible assoication for $\al_1'$ on the $\szfo$, so it is unique up to Dehn twists about $\al_1$.

\epsfysize=4 cm
\begin{figure}[htb]
\center{
\leavevmode
\epsfbox{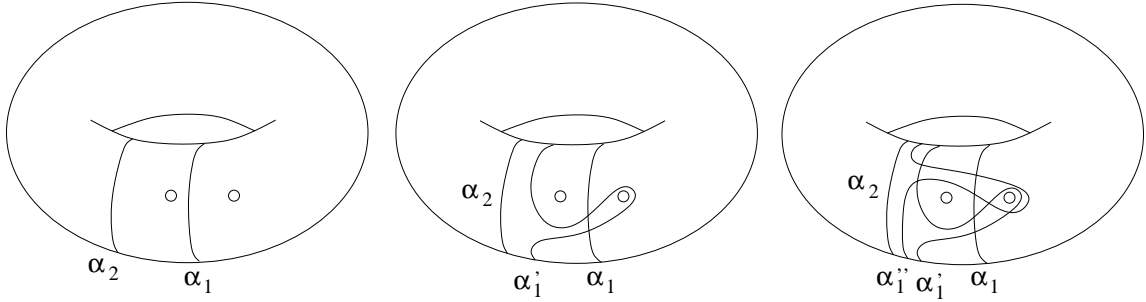}
\caption{Step 5. (i)  The fixed choices for $\al_1$ and $\al_2$ (ii) The unique choice for $\al_1'$.
(iii) One of the two choices for $\al_1''$.}           
}
\end{figure}

Since the vertices of $\cpo(\szfo)$ corresponding to $\al_1$, $\al_1'$ and $\al_1''$ are a quadrilateral
triple, there are two choices for $\al_1''$ (the argument is the same for both choices, so we treat one 
only).

Again, $\al_2'$ must have trivial intersection with $\al_1'$, and it must differ from $\al_1''$ and $\al_2$
by elementary moves.  The curves $\al_1'$ and $\al_1''$ divide $\sot$ into regions, two of which are
punctured disks, whose punctures come from the punctures of $\sot$.  Applying the Punctured Disks
Argument (see {\it Alternating Pentagons}), we see that there is no such curve $\al_2'$.  This is a
contradiction, so $\al_1'$ (and hence $\al_1''$) is separating on the $\sot$.



\step{6} We now have all of the information we need to finish the proof: $\al_1$ and $\al_2$ are
nonseparating, while $\al_1'$ and $\al_1''$ are separating on $\sot$.  Then, up to Dehn twists about
$\al_1$ and $\al_2$, $\al_1'$ and $\al_1''$ are the curves shown in the figure.  

\epsfysize=4 cm
\begin{figure}[htb]
\center{
\leavevmode
\epsfbox{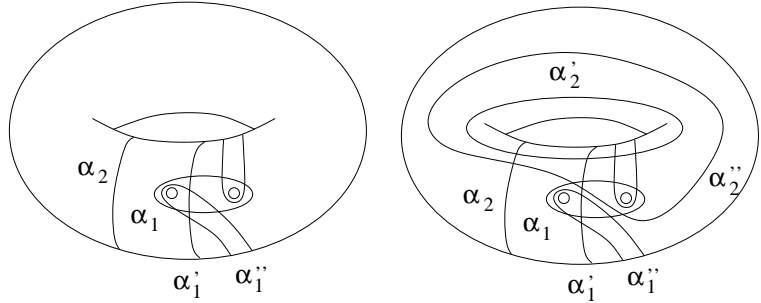} 
\caption{Step 6. (i) The unique choices for $\al_1$ and $\al_1''$. (ii) The unique choices for $\al_2'$
and $\al_2''$.}
}
\end{figure}

Finally, there are unique choices for $\al_2'$ and $\al_2''$, as $\al_2'$ must have trivial intersection
with $\al_1'$ and must have minimal intersection with $\al_2$ and $\al_1''$, while $\al_2''$ must have
trivial intersection with $\al_1''$ and must have minimal intersection with $\al_2$ and $\al_1'$. 

Note that the alternating pentagons described here correspond exactly to the pentagonal 2-cells in the
pants complex given by Figure~\ref{hexpic}.

Since triangles, alternating squares, alternating pentagons, and almost-alternating hexagons are all
preserved by automorphisms of the pants graph, and we have shown that these shapes correspond
precisely to the 2-cells of the pants complex, it follows that the automorphism groups of the
pants graph and the pants complex are the same.

\section{Proof of $\aut \cpos \simeq \mods$} \label{graph-curvecomplex}

For the general case of second isomorphism of the Main Theorem, we use a reciprocal strategy to that employed for the
5-times punctured sphere.  Instead of associating to a curve every pants decomposition containing that curve, we associate a collection of pants decompositions with exactly one of them containing the curve.

Recall that $S$ is a closed, oriented surface of negative Euler characteristic, which is not the twice punctured
torus.  

We are now ready to describe a map $\phi$ from $\aut\cpos$ to $\aut \ccz(S)$ (the
space of self-maps of the vertices of $\cc(S)$). We will then show that $\phi$ is well-defined, that $\phi$ naturally
extends to a map from $\aut\cpos$ to $\aut\cc(S)$, and finally that $\phi$ is an isomorphism between
$\aut\cpos$ and $\aut\cc(S)$.  Recall that we have already shown $\aut\cps \cong \aut\cpos$. 

Let $A$ be an automorphism of the pants graph of a surface $S$.  We define $\phi(A): \ccz(S) \ra
\ccz(S)$ (and hence $\phi$) by way of saying what $\phi(A)$ does to each vertex of $\cc(S)$: 

If $v$ is a vertex of $\cc(S)$ and $(F_v,X)$ is some marked Farey graph in the pants graph corresponding to $v$
(recall that there is a choice here; see Marked Farey Graphs, Section~\ref{graph-complex}), then $\phi(A)(v)$ is
defined to be
the unique vertex of the curve complex corresponding to the marked Farey graph $(A(F_v),A(X))$.

\paragraph{$\phi$ is well-defined.} Let $v$ be a vertex in the curve complex corresponding to the curve
$\al_1$.  Recall that the vertex $\phi(A)(v)$ is determined by the choice of marked Farey graph representing
$v$, and the choice of marked Farey graph comes down to a choice of pants decomposition containing $\al_1$. 
Therefore, we need to show that if $p$ and $p'$ are two pants decompositions which give rise to the marked
Farey graphs $(F_v,X)$ and $(F_v',X')$ representing $v$, then the curves corresponding to
$(A(F_v),A(X))$ and $(A(F_v'),A(X'))$ are the same. 

Actually, by the connectedness of the pants graph (of the surface obtained by cutting $S$ along
$\al_1$), we only need to treat the case when $p$ and $p'$ differ by a move, say $\al_2 \ra \al_2'$ (so
$X$ corresponds to $\{\al_1,\al_2,\dots,\al_n\}$ and $X'$ corresponds to
$\{\al_1,\al_2',\dots,\al_n\}$). 

\piece{Outline.} To insure that the vertices of the curve complex corresponding to $(A(F_v),A(X))$ and
$(A(F_v'),A(X'))$ are the same, we will find a way (intrinsic to the pants graph) of showing that the
vertices of $\cc(S)$ associated to ($F_v$,$X$)  and ($F_v'$,$X'$) are the same.  The idea is as follows: 
we will find a 2-curve small loop $\lp$ in $\cpos$, such that four of the vertices are (in order) 
$WXX'Y$, where $W$ and $X$ are in ($F_v$,$X$) and $X'$ and $Y$ are in ($F_v'$,$X'$).  Note that there is
no Farey graph in the pants graph containing either the triple of points $WXX'$ or the triple $XX'Y$
(since the pairs $WX'$ and $XY$ each only share $n-2$ curves).  In other words, $WXX'Y$ is an
alternating sequence of vertices in the pants graph.

\begin{figure}[htb]
\center{
\leavevmode
\epsfbox{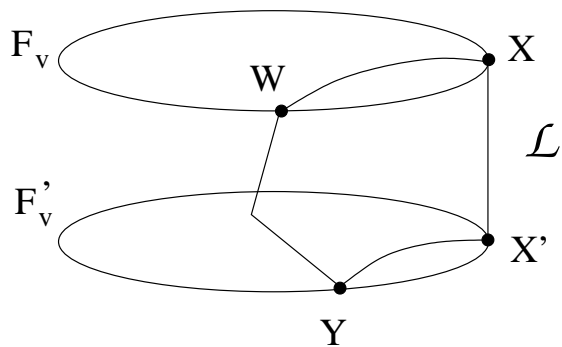} 
\caption{The 2-curve small loop.}
}
\end{figure}

Since the property of alternating was defined completely in terms of the geometry of the pants graph, it
must be preserved by the automorphism $A$, so the sequence $A(W)A(X)A(X')A(Y)$ must also be alternating.
Say the edge $A(W)A(X)$ corresponds to an elementary move $\be_1' \ra \be_1$ (say
$\{\be_1',\be_2,\dots,\be_n\} \ra \{\be_1,\dots,\be_n\}$)  for some curves $\be_1'$ and
$\be_1$---so the marked Farey graph $(A(F_v),A(X))$ is associated to $\be_1$, then (by the alternating
property) the edge $A(X)A(X')$ must correspond to an elementary move $\be_2 \ra \be_2'$, where $\be_2$ and
$\be_2'$ are some other curves (so $A(X')$ corresponds to $\{\be_1,\be_2',\dots,\be_n\}$. Now,
the edge $A(X')A(Y)$ cannot correspond to a move $\be_2' \ra \star$ (by the alternating property again),
nor can it correspond to any other move $\be_i \ra \star$, for $i > 2$, since $\lp$ was chosen to be a
2-curve small loop, and in such a loop, all but two of the curves stay fixed. Therefore, the edge
$A(X')A(Y)$ corresponds to a move $\be_1 \ra \star$, and so the marked Farey graph $(A(F_v'),A(X'))$
corresponds to $\be_1$, just like $(A(F_v),A(X))$. 

\piece{Illegal Moves.} Suppose $p$ and $p'$ differ by the move $\al_2 \ra \al_2'$.  There is one {\em illegal move}
which we will rule out now---the case when $\al_1$ and $\al_2$ lie on a $\sot$, and $\al_1$, $\al_2$, and $\al_2'$ are
all non-separating on the $\sot$. The reason we need to rule this out is that this type of move does not appear in the
almost-alternating hexagon (or any of the other small loops previously discussed), and this makes it difficult to find
the loop $\lp$ in this case. We need not concern ourselves with this case, because the pants graph is still connected
even when the edges corresponding to these moves are removed.  In other words, there is a pants decomposition $p''$
which differs by a move from both $p$ and $p'$ so that the above {\em illegal move} does not occur in passing from
$p$
to $p''$ or from $p''$ to $p'$. We now prove this fact:

Since $\al_1$ is nonseparating on $\sot$, its complement in the subsurface $\sot$ is a
$\szfo$, with the four ``punctures'' corresponding to the two punctures from $\sot$ and the two ``sides''
of $\al_1$.  Since $\al_2$ and $\al_2'$ are nonseparating curves on the $\szfo$ which differ by
an elementary move, they each associate (see Section~\ref{complexes}) different punctures with a given side of $\al_1$
(in particular, they don't associate the two punctures of the $\sot$).  Since any edge in $\cpo(\szfo)$ is part of a
triangle, there is a vertex (two in fact) in $\cpo(\szfo)$ (corresponding to a curve $\al_2''$) connected to both
vertices corresponding to $\al_2$ and $\al_2'$, where $\al_2''$ associates the two punctures of $\sot$ (i.e. is
separating on $\sot$), since the pants decompositions corresponding to the three vertices of a triangle of
$\cpo(\szfo)$ have all three different associations. Choose $p''$ so that the moves $p \ra p''$ and $p'' \ra p$
correspond to the moves $\al_2 \ra \al_2''$ and $\al_2'' \ra \al_2'$.

\epsfysize=4 cm
\begin{figure}[htb]
\center{
\leavevmode
\epsfbox{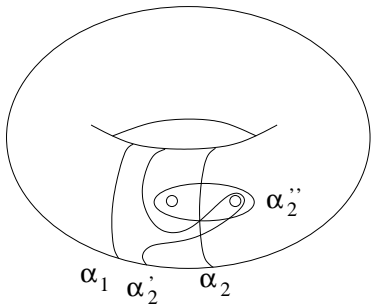} 
\caption{Circumventing illegal moves.}
}
\end{figure}

\piece{Finding the Loop.} In order find the loop $\lp$, there are three cases to consider:

1. $\al_1$ and $\al_2$ lie on disjoint subsurfaces

2. $\al_1$ and $\al_2$ lie on a $\szfi$

3. $\al_1$ and $\al_2$ lie on a $\sot$

Case 1: Choose $\lp$ to be the alternating square $WXX'Y$, where $W$, $X$, $X'$, and $Y$ have the same
curves (in name) in their corresponding pants decompositions as $P$, $Q$, $R$, and $S$ (from the
topological description of alternating squares).

Case 2: Choose $\lp$ to be the alternating pentagon $WXX'YZ$ (for unique $Z$), where $W$, $X$, $X'$,
$Y$, and $Z$ have the same curves in their corresponding pants decompositions as $P$, $Q$, $R$,
$S$, and $T$ respectively (from the topological description of alternating pentagons). 

Case 3: Since we eliminated illegal moves, there are two possibilities for the curves $\al_1$,
$\al_2$, and $\al_2'$---either (i) $\al_1$ and $\al_2$ is nonseparating, and $\al_2'$ is separating, or (ii)
$\al_1$ is separating, $\al_2$ and $\al_2'$ are nonseparating.  By symmetry, the first possibility takes
care of the case when $\al_1$ is nonseparating, $\al_2$ is separating, and $\al_2'$ is nonseparating. 

For the first possibility, choose $\lp$ to be the almost alternating hexagon $WXX'YZA$ (for unique $Z$
and $A$) where $W$, $X$, $X'$, $Y$, $Z$, and $A$, correspond to $R$, $S$, $T$, $U$, $P$, and $Q$,
respectively (from the topological characterization of almost alternating hexagons).  For the latter
case, $W$, $X$, $X'$, $Y$, $Z$, and $A$, should correspond to $S$, $T$, $U$, $P$, $Q$, and $R$,
respectively. 

\paragraph{$\phi(A)$ is an automorphism.} Since the curve complex has the property that every set of $k+1$
mutually connected vertices is the $1$-skeleton of a $k$-simplex in the curve complex, we only need to check
that $\phi(A)$ extends to an automorphism of the $1$-skeleton of the curve complex, i.e. that $\phi(A)$
takes vertices connected by edges to vertices connected by edges.

Suppose that $v$ and $w$ are vertices of the curve complex corresponding to curves $\al$ and $\be$ on
the surface. We will show that $v$ and $w$ are connected by an edge in the curve complex if and only if
there are marked Farey graphs $(F_v,P)$ and $(F_w,P)$ corresponding to $v$ and $w$, where the Farey graphs
share a unique point $P$, which is also the marked point on both Farey graphs. Since the number of
intersection points between Farey graphs is preserved under the pants graph automorphism $A$, it follows
that edges are preserved under $\phi(A)$. 

Let $P$ be any pants decomposition containing $\al$ and $\be$, say $\{\al,\be,\ga_1,\dots,\ga_n\}$.  Then
let $F_v$ and $F_w$ be the set of all pants decompositions of the form $\{\star,\be,\ga_1,\dots,\ga_n\}$ and
$\{\al,\star,\ga_1,\dots,\ga_n\}$.  Then $(F_v,P)$ and $(F_w,P)$ have the desired property.  This
construction is possible if and only if $\al$ and $\be$ are part of some pants decomposition, which is
equivalent to $v$ and $w$ being connected by an edge.

\paragraph{$\phi$ is an isomorphism.}  We show that $\phi$ is a homomorphism that is injective and
surjective.

\bigskip

\noindent {\em Homomorphism.} We will show that $\phi(AB) = \phi(A)\phi(B)$ in $\aut\cc(S)$ by showing that
$\phi(AB)v = \phi(A)\phi(B)v$ for any vertex $v$ in $\cc(S)$.  $\phi(AB)v$ is, by definition, the vertex in $\cc(S)$
corresponding to $(AB(F_v),AB(P))$, where $(F_v,P)$ is a marked Farey graph in $\cps$ representing $v$, and $AB$,
being an element of $\aut\cps$ acts on it in the natural way.  On the other hand, $\phi(B)v$ is the vertex $w$ in
$\cc(S)$ corresponding to $(B(F_v),B(P))$, and hence $\phi(A)\phi(B)v$ is the vertex of $\cc(S)$ corresponding to
$(A(F_w),A(Q))$, where $(F_w,Q)$ is some Farey graph representing $w$.  We can in fact choose $(F_w,Q)$ to be
$(B(F_v),B(P))$, and so $\phi(A)\phi(B)v$ is $(AB(F_v),AB(P))$, which is the same as $\phi(AB)v$.

\bigskip

\noindent {\em Injectivity.} Suppose $\phi(A)$ is the identity automorphism of the curve complex, and let $P$ be a
vertex of the pants graph corresponding to the pants decomposition $\{\al_1,\dots,\al_n\}$, where
$v_1,\dots,v_n$ are the vertices of the curve complex corresponding to the $\al_i$. We will show that
$A(P) = P$.  There is a unique choice of $F_{v_i}$ (the Farey graph corresponding to the pants
decompostions $\{\al_1,\dots,\al_{i-1},\star,\al_{i+1},\dots,\al_n\}$) representing the $v_i$ so that
the $F_{v_i}$ all intersect at the point $P$ in the pants graph. 

Since $A(F_{v_1}),\dots,A(F_{v_n})$ must be marked Farey graphs corresponding to $v_1,\dots,v_n$, it follows
that their common intersection point is $P$.  Thus $A(P) = P$.

\bigskip

\noindent {\em Surjectivity}: There is a natural map $\psi : \aut\cc(S) \lra \aut\cps$ obtained by using the
isomorphism $\eta: \aut\cc(S) \lra \mcg(S)$ of Theorem A.  The isomorphism $\eta$ is natural in the
sense that if a curve $\al$ corresponds to a vertex $v$ of the curve complex, then the curve
$\eta(\mathcal{A})(\al)$ represents the vertex $\mathcal{A}(v)$ for any $\mathcal{A} \in \aut\cc(S)$.

Let $\mathcal{A} \in \aut\cc(S)$, and let $P$ be a vertex of the pants graph, where $P$ is the set of curves
$\{\al_1,\dots,\al_n\}$.  Then $\psi(\mathcal{A})(P)$ is defined to be the
vertex of the pants graph which is the curves
$\{\eta(\mathcal{A})(\al_1),\dots,\eta(\mathcal{A})(\al_n)\}$.  Note that $\psi(\mathcal{A})$ preserves
edges in $\cps$ since it is has the same action as a mapping class element.

It remains to show that $\phi \circ \psi (\mathcal{A}) = \mathcal{A}$, where $\mathcal{A} \in \aut\cc(S)$. If $v$ is a
vertex of the curve complex, then we need to show that $\psi(\mathcal{A})(F_v,P) = (F_{\mathcal{A}(v)},Q)$, where
$(F_v,P)$ and $(F_{\mathcal{A}(v)},Q)$ are marked Farey graphs in the pants graph corresponding to the vertices $v$
and $\mathcal{A}(v)$ of the curve complex.  If $v$ corresponds to a curve $\al$, then the vertices of $F_v$ are all
the pants decompositions of the form $\{\star,\al_2,\dots,\al_n\}$ for some curves $\al_2,\dots,\al_n$ which have
trivial intersection with $\al$.  By definition, the vertices of $\psi(\mathcal{A})(F_v)$ are all of the pants
decompositions of the form $\{\star,\eta(\mathcal{A})(\al_2),\dots,\eta(\mathcal{A})(\al_n)\}$, and
$\psi(\mathcal{A})(P)$ is the pants decomposition
$\{\eta(\mathcal{A})(\al),\eta(\mathcal{A})(\al_2),\dots,\eta(\mathcal{A})(\al_n)\}$. Since $\eta(\mathcal{A})(\al)$
corresponds to the vertex $\mathcal{A}(v)$ of the curve complex, we are done.

\bibliography{pants}

\bibliographystyle{plain}

\noindent
Dan Margalit\\
Dept. of Mathematics, University of Chicago\\
5734 University Ave.\\
Chicago, Il 60637\\
E-mail: juggler@math.uchicago.edu

\end{document}